\documentclass[11pt]{article}
\usepackage[latin1]{inputenc}
\usepackage{tikz}
\usepackage{amssymb}
\usepackage{amsmath}
\usepackage{amsthm}

\newcommand{\R}{\mathbb{R}}
\newcommand{\Q}{\mathbb Q}

\newcommand{\Z}{\mathbb Z}

\theoremstyle{definition}

\newtheorem{thm}{Theorem}[section]
\newtheorem{prop}[thm]{Proposition}
\newtheorem{lem}[thm]{Lemma}
\newtheorem{rem}[thm]{Remark}
\newtheorem{defi}[thm]{Definition}
\newtheorem{cor}[thm]{Corollary}

\newcommand{\vs}{\vspace{0.3cm}}

\newcommand{\vsp}{\vspace{0.1cm}}

\date{}
\author{}

\begin{document}

\title{On the space of left-orderings of virtually solvable groups}
\author{Cristobal Rivas\footnote{The first author was supported by CONICYT and labex
MILYON.} \& Romain Tessera\footnote{The second author was supported
by ANR BLAN AGORA  and ANR BLAN GGAA.}} \maketitle
\begin{abstract}
We show that the space of left-orderings of a countable virtually
solvable group is either finite or homeomorphic to a Cantor set. We
also provide an explicit description of the space of left-orderings of
$SOL=\Z^2\rtimes_T\Z$.
\end{abstract}

\section{Introduction}

The space of left-orderings $\mathcal{LO}(G)$ of a left-orderable
group $G$ is the set of all possible left-orderings on $G$ endowed
with a natural topology that makes it compact, Hausdorff and totally
disconnected, see \cite{sikora} or \S \ref{sec LO}. It was proved by
Linnell that this space is either finite or uncountable
\cite{linnell}. The problem of relating the topology of
$\mathcal{LO}(G)$ with the algebraic structure of $G$ has been of
increasing interest since the discovery by Dubrovina and Dubrovin
that the space of left-orderings of the braid groups is infinite and
yet contains isolated points \cite{dd}. Recently, more examples of
groups showing these two behaviors have appeared in the literature
\cite{dehornoy, ito dehornoy like, ito constructions, navas-hecke}.
All these groups admit a non-abelian free subgroup. However,
non-abelian free groups \cite{McCleary}, and more generally
non-trivial free products of groups have no isolated left-orderings
\cite{rivas free prod}.

In any case, the class of groups having isolated orderings is far
from being well understood. By contrast, left-orderable groups
admitting only finitely many left-orderings have been classified by
Tararin \cite[Theorem 5.2.1]{kopytov}. For short, we shall call
these groups ``Tararin groups". Tararin groups form a very
restrictive and easy to describe class of
finite-rank-solvable\footnote{A solvable group has finite rank, if
in its derived series all successive quotients have finite rank as
Abelian groups, see \cite{robinson}.} groups, see \S \ref{sec
tararin}. As already mentioned, these are the only known examples of
left-orderable {\it amenable} groups with isolated points. Could
they be the only ones? This paper brings in a modest contribution to
this problem, answering it positively for (virtually) solvable
groups.

\begin{thm}\label{teo main} {\em The space of left-orderings of
a countable virtually solvable group is either finite or a Cantor
set.}
\end{thm}

In particular we deduce
\begin{cor} {\em If $\Gamma$ is a (countable) left-orderable virtually solvable
group of infinite rank, then $\mathcal{LO}(\Gamma)$ is a Cantor
set.}
\end{cor}

Among Tararin groups, those which are virtually
polycyclic\footnote{A polycyclic group, is a solvable group whose
successive Abelian quotients, in its derived series, are finitely
generated, see \cite{raghunathan}.} turn out to be virtually
nilpotent, therefore we deduce

\begin{cor} {\em If $\Gamma$ is a left-orderable virtually polycyclic group of exponential
growth, then $\mathcal{LO}(\Gamma)$ is a Cantor set.}
\end{cor}

The study of left-orderable amenable groups is intimately related to
that of Conradian orderings (see \S \ref{sec conrad} for a
definition) and to local indicability. Recall that Morris gave a
beautiful proof of the fact that all left-orderable amenable groups
admit Conradian orderings  \cite{Mo} (see also  \cite{MR,ChKrop} for
older results in that direction, and \cite{Der} for an interesting
alternative proof).  Together with a fundamental observation of
Conrad \cite{conrad}, this provides a very natural characterization
of left-orderable amenable groups as those which are locally
indicable,  i.e.\ all their finitely generated subgroups have a
non-trivial morphism to $\Z$.

\vsp

The dichotomy shown in Theorem \ref{teo main} reminds of a similar
one, this time for all groups but in restriction to Conradian
orderings. Indeed, in
\cite{rivas: conrad} the first author proved that the space of
Conradian orderings of a countable group is either finite or
homeomorphic to the Cantor set. This implies for instance, the
general dichotomy for left-orderable groups having only Conradian
orderings, such as groups of sub-exponential growth \cite{Navas}.

\subsection*{Remarks about the proof}

It is now well known that countable groups are left-orderable if and
only if they act faithfully by order-preserving homeomorphisms on
the real line (see for instance \cite{ghys} or \S \ref{sec dynamical
real}). This dynamical approach have shown to be fruitful when
trying to understand the topology of the space of left-orderings of
a given group or family of groups, see for instance
\cite{Navas,rivas: conrad,rivas free prod}. In this work, this point
of view will be crucial.

\vsp

In a previous version of this paper, Theorem \ref{teo main} was
proved for the class of finite-rank solvable groups. An important
feature of this subclass is that they are virtually
nilpotent-by-Abelian \cite{robinson}. In this case, the proof can be
summarized as follows: if the ordering is not Conradian (in which
case, the conclusion follows from \cite{Navas}), then one can prove
that the ordering is induced, up to semiconjugacy, from an affine
action on the real line. It is then not to hard to see that the
underlying ordering is non-isolated, see Proposition \ref{prop
affine} and Corollary \ref{cor:affine}. This approach strongly
relies on the work of Plante \cite{plante}, who shows that any
action of a finite-rank solvable group on the real line quasi
preserves a Radon measure. However, as already noted by Plante,
there are actions of solvable groups on the real line (such as
$\Z\wr\Z$) in which no non-trivial Radon measure is quasi-preserved.
This situation is much more subtle and requires a careful dynamical
analysis. Roughly speaking, the proof consists in showing that the
action behaves ``at a certain scale" like an affine action (see
Section \ref{sec Plante} and Lemma \ref{lem no strong crossings}).

\subsection*{Organization of the paper}

In \S \ref{preliminars} we give some necessary background. Most of
the material is well known, but not all. Notably Corollary \ref{cor
no strong crossing} plays a crucial role in our proof. In \S
\ref{sec examples} we prove Theorem \ref{teo main} in a simple
example\footnote{Solvable Baumslag-Solitar groups could have been
treated in a similar way, we leave the easy adaptation of the proof
to the reader.}, namely $SOL$. We also give an explicit description
of its space of left-orderings. In \S \ref{sec Plante} we illustrate
the difficulties arising when dealing with solvable groups of
infinite rank on a specific example, due to Plante, of an action of
$\Z\wr \Z$ on the line. Finally the proof of Theorem \ref{teo main}
is carried out in \S \ref{sec Proof of Main thm}.
\begin{rem}
In various places we employ the useful terminology
``pseudo-ordering" (Definition \ref{pseudo}) which means an
invariant ordering on a quotient of $G$ by a (not necessarily
normal) subgroup, but seen as a partial ordering on $G$.
\end{rem}
\vs

\noindent {\bf Acknowledgements:} The first-named author would like
to thanks the hospitality of ENS-Lyon where he enjoys a one year and
a half posdoctoral position. He is specially grateful to M.
Triestino, A. Glutsyuk and \'E. Ghys for interesting discussions
around this and other related problems.


\section{Preliminaries}\label{preliminars}


\subsection{The topology on $\mathcal{LO}(G)$}

\label{sec LO}

A basis of neighborhoods in $ \mathcal{LO} (G)$ is the family of the
sets $\, V_{f_1,\ldots,f_k}:=\{\preceq \mid id\preceq f_1, \ldots ,
id\preceq f_k\} \,$, where $\{f_1,\ldots,f_k\}$ runs over all finite
subsets of $G$. If $G$ is countable, then this topology is
metrizable. For instance, if $G$ is finitely generated, we may
define $\,dist(\preceq,\preceq') = 1 / 2^n$, where $n$ is the first
integer such that $\, \preceq \,$ and $\, \preceq'\,$ do not
coincide on $n$-th ball (with respect to some generator system). For
each $g\in G$ and $\preceq\in \mathcal{LO}(G)$, one can define an
other element $g(\preceq)\in \mathcal{LO}(G)$ whose {\em positive
cone} is the set of elements $f\in G$ such that $gfg^{-1}\succ id$.
This defines a continuous representation of $G$ in
$Homeo(\mathcal{LO}(G))$ called action by {\it conjugation} of $G$
on its space of left-orderings.

The following definition is classical \cite{botto
rhemtulla,kopytov}. Given a left-ordered group $(G,\preceq)$ and a
subgroup $H$, we say that $H$ is {\em convex} if for every $g\in G$
such that $h_1\preceq g\preceq h_2$, for $h_1, \; h_2$ in $H$, we
have that $g\in H$. Convex subgroups have the nice property that
they induce a total ordering on the left-cosets $G/H$ by $$g_1H
\preceq^* g_2H \Leftrightarrow g_1h_1 \preceq g_2 h_2 \; \text{ for
all } h_1,h_2 \text{ in } H.$$ This ordering is invariant under the
$G$-action by left translation on $G/H$. In particular, if $H$ is
normal, then $G/H$ is a left-orderable group, in which case we call
$\preceq^*$ the {\em projected} or {\em quotient} ordering. It
follows that $\preceq$ decomposes ``lexicographically" as the
ordering on the $H$-cosets and the ordering restricted to $H$. More
precisely we have that $$ id \prec g \Leftrightarrow \left\{
\begin{array}{cc} H\prec^* gH, \text{ or}  \\  H=gH \text{ and }
id \prec g
\end{array}\right.$$

Elaborating on this, we conclude (see \cite{rivas dichotomy} for
more details)

\begin{prop} \label{prop convex}{\em Let $\preceq$ be a
left-ordering on $G$ and let $H$ be a convex subgroup. Then there is
a continuous injection $\mathcal{LO}(H)\to \mathcal{LO}(G)$, having
$\preceq$ in its image. Moreover, if in addition $H$ is normal, then
there is a continuous injection $\mathcal{LO}(H)\times
\mathcal{LO}(G/H)\to \mathcal{LO}(G)$ having $\preceq$ in its
image.}
\end{prop}

\begin{rem}\label{rem:ConradExtension} Let $\preceq$ be a left-ordering
on $G$, and $H$ a normal convex subgroup. Then, it is not hard to
check that if the restriction of $\preceq$ to $H$ and the projection
of $\preceq$ to $G/H$ are Conradian, then $\preceq$ is also
Conradian.
\end{rem}

Observe that the set of convex subgroups of a given left-ordering,
are totally ordered for the inclusion. We call it the convex series
of $G$. The following Corollary is well known, and appears for
instance in \cite{rivas: conrad}. For the readers' convenience we
sketch a proof. Let $G$ be a group and $H$ be a (not necessarily
normal) subgroup. Then to any $G$-invariant ordering $\preceq$ on
$G/H$, one can define the so-called {\em opposite} ordering
$\preceq^{op}$ defined by $g\prec^{op} f \Leftrightarrow g\succ f$.
We sometimes say that $\preceq^{op}$ is obtained by {\em flipping}
$\preceq$. Clearly $\preceq^{op}$ is also $G$-invariant.

\begin{cor} \label{cor infinite convex series}{\em Suppose that a left-ordered
group $(G,\preceq)$ has infinitely many convex subgroup. Then
$\preceq $ is non-isolated.}
\end{cor}

\noindent{\em Sketch of the proof:} If the convex series is
infinite, then either there exists an infinite increasing sequence
of convex subgroups $C_1<C_2\ldots$ or an infinite decreasing
sequence $C'_1>C'_2\ldots$. Flipping the ordering on $C_{n+1}/C_n$
(resp. on $C'_n$), one obtains a sequence of orderings $\preceq_n$,
distinct from $\preceq$, which converges to $\preceq$ when $n$ goes
to infinity. $\hfill\square$

\vsp

Note that we often have to deal with orderings on a quotient $G/H$,
where $H$ is not necessarily normal. It will therefore be convenient
to see $G$-invariant orderings on $G/H$ as a ``pseudo-ordering" on
$G$, namely
\begin{defi}\label{pseudo} (Pseudo-ordering) A pseudo-ordering $\preceq$ on a group $G$ is a
left-invariant partial ordering on $G$ induced from a $G$-invariant
ordering $\preceq^*$ on a quotient $G/H$: $id\prec g\Leftrightarrow
H\prec^* gH$. Given a pseudo-ordering on $G$, the set of elements
which are not comparable to the neutral element coincides with the
subgroup $H$.
\begin{itemize}
\item (Convex subgroup) A convex subgroup $C$ of a pseudo-ordered
group is defined similarly as for orderings, with the additional
requirement that $C$ must contain $H$. Note that  $H$ itself is
convex, and so is the minimal convex subgroup.

\item (Quotient by a convex subgroup) Given a pseudo-ordering
$\preceq$ on $G$ and $C$ a convex subgroup, the pseudo-ordering
on $G$ induced from the ordering on $G/C$ will be called the
quotient of $\preceq$ by $C$.
\end{itemize}
\end{defi}

\begin{rem}\label{rem convex}
If $\preceq'$ is the quotient of $\preceq$ by a convex subgroup $C$
of $(G,\preceq)$, then convex subgroups of $(G,\preceq')$ are
exactly those convex subgroups of $(G,\preceq)$ containing $C$.
\end{rem}
\subsection{Tararin groups}
\label{sec tararin}

We give a slight modification of the original statement of Tararin
\cite[Theorem 5.2.1]{kopytov}, describing groups admitting only
finitely many left-orderings. Recall that a series
$$\{1\}=G_0\lhd G_1 \lhd \ldots \lhd G_m=G,$$ is said to be
{\em rational} if each quotient $G_{i+1}/G_i$ is torsion-free
rank-one Abelian.

\begin{thm}[Tararin] \label{teo tararin}{\em Let $G$ be a left-orderable
group. If $G$ admits only finitely many left-orderings, then $G$
admits a unique (hence characteristic) rational series
$$\{1\}=G_0\lhd G_1 \lhd \ldots \lhd G_m=G,$$ such that, for every
$2\leq i\leq m$, there is an element of $G_i/G_{i-1}$ whose action
by conjugation on $G_{i-1}/G_{i-2}$ is by multiplication by a
negative rational number. We shall call such a group a Tararin
group. }
\end{thm}

\begin{rem} \label{rem orderes on Tararin's}
The left-orderings on a Tararin group $G$ are very easy to describe.
Indeed, if $\{1\}=G_0\lhd \ldots \lhd G_m=G,$ is the associated
rational series, then on each quotient $G_i/G_{i-1}$, being rank-one
torsion free Abelian, there is --up to flipping-- a unique
left-ordering coming from an embedding into $\Q$. For every $i$, let
$\preceq_i$ be a choice of an ordering on $G_i/G_{i-1}$. Then we can
produce a left-orderings on $G$ by declaring
$$g\succ id \Leftrightarrow \left\{\begin{array}{lc} gG_{m-1}\succ_m
G_{m-1}\;, \text{ or }\\ g\in G_{m-1}\; , \text{ and }\;
g\succ_{m-1}
id\, , \text{ or }\\ \hspace{1.5 cm}\vdots\\
g\in G_1\; , \text{ and }\; g\succ_1 id.
\end{array}\right.$$
It is not hard to check that in this way we can produce all possible ($2^m$)
left-orderings (in fact, it is easy to show that they are all
Conradian). Moreover, in any such ordering, the groups $G_i$ are
convex, and conversely, every convex subgroup is of this form.
\end{rem}

\begin{cor} {\em Let $G$ be a virtually polycyclic group admitting only
finitely many left-orderings. Then it admits a unique filtration
such that $G_i/G_{i-1}\simeq \Z$. The action of (the generator of)
$G_i/G_{i-1}$ on $G_{i-1}/G_{i-2}$ is by multiplication by $-1$.}

\end{cor}

Since in a virtually polycyclic group, the group generated by
$\{g^2\mid g\in G\}$ has finite index \cite{raghunathan}, we
deduce

\begin{cor}\label{coro Tararin finite-rank-solvable}{\em A virtually
polycyclic group having only finitely many left-orderings is
virtually nilpotent.}
\end{cor}

We finish this section with a rigidity statement for actions of
Tararin groups on the line. Let $G$ be a Tararin group and
$\{1\}=G_0\lhd G_1 \lhd \ldots \lhd G_m=G$ its associated rational
series. Since $G_{m-1}$ is convex in every left-ordering of $G$, and
$G/G_{m-1}$ is Abelian we have that the sign of $\gamma\in
G\setminus G_{m-1}$ is preserved under conjugation. On the other
hand, Remark \ref{rem orderes on Tararin's} says that starting from
a left-ordering $\preceq$, any other left-ordering on $G$ is
obtained by flipping the ordering on some of the convex subgroups.
It follows from Theorem \ref{teo tararin}, that any flipping on any
proper convex subgroup can be realized as conjugations by some
element in $G$. This shows

\begin{prop}\label{prop two orbits of Tararin}
{\em The conjugation action of a Tararin group $G$ on
$\mathcal{LO}(G)$ has two orbits. Moreover, for any pair of
left-orderings $\preceq$, $\preceq^\prime$ of $G$ there is $g\in G$
such that $g(\preceq)$ and $\preceq^\prime$ coincide over $G_{m-1}$.
Moreover, if we let $\gamma_T$ be an element in $G\setminus G_{m-1}$
which acts on $G_{m-1}/G_{m-2}$ as multiplication by a negative
number, then $g$ can be taken either in $G_{m-1}$ or in
$\gamma_TG_{m-1}.$}
\end{prop}

\subsection{Dynamical realization}

\label{sec dynamical real}

As mentioned in the introduction, an important ingredient in our proof
of Theorem \ref{teo main} is the fact that countable left-orderable
groups naturally act by order-preserving automorphism of the real
line, and vice versa, a group acting faithfully by order-preserving
automorphism of the real line is left-orderable \cite{ghys}.

More precisely, given a left-ordered group $(G,\preceq)$, there is
an embedding of $G$ into $Homeo_+(\R)$, the group of order
preserving automorphism of the real line, such that: \begin{itemize}

\item $G$ acts without global fixed points,

\item for $f,g$ in $G$, we have that $f\prec g \Leftrightarrow
f(0)< g(0)$, and

\item the set of fixed points of a non-trivial $f\in G$ has empty
interior.

\end{itemize}
This construction extends to pseudo-orderings in the sense that
every pseudo-order on $G$ with minimal convex subgroup $H$ can be
induced from an action of $G$ on the real line, where $H$ is the
stabilizer of $0$.

We call such an action, a {\em dynamical realization} of
$(G,\preceq)$. Conversely, given an embedding of $G$ into
$Homeo_+(\R)$, we can induce a left-ordering as follows: take
$(x_1,x_2,\ldots )$ a dense sequence in $\R$, and declare that an
element  $g\succ_{(x_1,x_2,\ldots)} id$ if and only if $g(x_i)>x_i$,
where $i$ is such that $g(x_j)=x_j$ for every $j<i$. We call such an
ordering, an {\em induced ordering} from the action. Note that with
this procedure we can recover a left-ordering from its dynamical
realization by taking $x_1=0$.

\begin{rem}\label{rem covariant actions}
Let $g\in G$, and $\preceq=\preceq_{(x_1,x_2,\ldots)}$ be an
ordering  induced from an action of $G$ on the real line. Then
$g(\preceq)$ is the left-ordering induced from the sequence
$(g^{-1}(x_1),g^{-1}(x_2),\ldots)$.

\end{rem}

\subsection{Orderings induced by affine actions}

A general procedure for trying to approximate a given left-ordering
$\preceq$ on a countable group, is to consider its dynamical
realization, and to induce an ordering $\preceq^\prime$ from a
sequence $(x_1,x_2,\ldots)$ where $x_1$ is close to $0$. The fact
that $\preceq^\prime$ is close to $\preceq$ when $x_1$ is close to
$0$ follows from the continuity of the action, and from the fact
that $0$ has a free orbit (details are left to the reader). The
problem however, is that the two orbits may induce the same ordering. This
is the case for instance if $\preceq$ has a least positive element.

Our first step in proving Theorem \ref{teo main} is that
left-orderings induced from non-Abelian affine actions, are
non-isolated.

\begin{prop}\label{prop affine} {\em Let $\Gamma$ be a countable group, and
suppose a left-ordering $\preceq $ on $\Gamma$ is induced from a
faithful (order-preserving) affine action on the real line. Then, if
$\Gamma$ is non-Abelian, $\preceq$ is approximated by its
conjugates.}
\end{prop}

\noindent{\em Proof:} Let $\preceq=\preceq_{(x_1,x_2,\ldots)}$ be
the left-ordering induced from the sequence $(x_1,x_2,\ldots)$. We
note that, since the elements in the affine group have at most one
fixed point, it is enough to specify two points. So
$\preceq=\preceq_{(x_1,x_2)}$.

By assumption, $\Gamma$ has both non-trivial homotheties and
non-trivial translations. It follows that the subgroup made of
translations has dense orbits. In particular, the countable set
$\Omega$ consisting of the points in $\R$ which are fixed by some
non-trivial element (homothety) of $\Gamma$, is also dense in $\R$.
Therefore, given any two points $x_1,y_1\in \R$, $x_1\not=y_1$,
there is a non-trivial homothety $h\in \Gamma$, having its unique
fixed point between $x_1$ and $y_1$. In particular, the
left-orderings $\preceq_{(x_1,x_2)}$ and $\preceq_{(y_1,y_2)}$,
induced from $(x_1,x_2)$ and $(y_1,y_2)$ respectively, are
different.

We now show that $y_1$ may be chosen so that $\preceq_{(y_1,y_2)}$
is close to $\preceq_{(x_1,x_2)}$. As noted earlier, this is obvious
if $x_1$ has a free orbit, so let us suppose that it is not the
case. Let $Stab_\Gamma(x_1)$ be the stabilizer of $x_1$ in $\Gamma$,
and let $S\subset \Gamma$ be a finite set of $\preceq$-positive
elements. We write $S=S_1\cup S_2$ where $S_1=S\cap
Stab_\Gamma(x_1)$, and we call $I$ the open interval between $x_1$
and $x_2$. Since $S_2$ is finite, there is a small neighborhood $U$
of $x_1$ such that $\gamma(x)>x$ for every $x\in U$ and every
$\gamma\in S_2$. On the other hand, for every $\gamma \in S_1$, we
have that $\gamma(x)>x$ for every $x\in I$ (recall that $\gamma\in
S_1$ is an homothety). Thus, if we take $y_1\in I\cap U$ then
$\preceq_{(y_1,y_2)}$ and $\preceq$ ($y_2$ being any point) are in
the same open set associated to $S$. Since $I\cap U$ has non-empty
interior, it is easy to see that $(y_1,y_2)$ may be chosen so that
$(y_1,y_2)=(\gamma(x_1),\gamma(x_2))$ for some $\gamma\in \Gamma$,
which shows that $\preceq$ is approximated by its conjugates.
$\hfill\square$

\vsp

To state the following corollary, recall that given two actions
$A_1$, $A_2$ of a group $\Gamma$ on the real line, we say that $A_1$
is {\em semi-conjugated} to $A_2$ if there is an increasing
surjective function $F:\R\to\R$ such that
$$F\circ A_1(\gamma)=A_2(\gamma)\circ F \;\; (\forall \gamma\in \Gamma).$$

\begin{cor}\label{cor:affine}
{\em Let $(\Gamma,\preceq)$ be a countable, left-ordered group.
Suppose there is an affine (order-preserving) action $A:\Gamma\to
Aff_+(\R)$ whose kernel is convex in $\preceq$, and whose range is
non-Abelian. Suppose further that the dynamical realization of
$(\Gamma,\preceq)$ is semi-conjugated to $A$. Then $\preceq$ is
non-isolated.}
\end{cor}

\noindent {\em Proof:} In light of Proposition \ref{prop convex}, it
is enough to show the corollary when $A$ is a faithful action. Let
$F:\R\to\R$ be the function that realizes the semi-conjugation, and
let $Stab_{A(\Gamma)}(F(0))$ be the stabilizer of $F(0)$ in
$A(\Gamma)$: this is an Abelian subgroup of $\Gamma$. We claim that
it is convex. Indeed, if $\gamma_1\prec g\prec \gamma_2$, then
$\gamma_1(0)<g(0)<\gamma_2(0)$. So $F(\gamma_1(0))\leq F(g(0))\leq
F(\gamma_2(0))$, and hence $A(\gamma_1)(F(0))\leq A(g)(F(0))\leq
A(\gamma_2)(F(0))$. From where the claim follows.

Now, if $Stab_{A(\Gamma)}(F(0))$ is trivial then Proposition
\ref{prop affine} applies directly, since in this case $\preceq$ is
realized as the induced ordering from $F(0)$ in the action $A$. If
it has rank $>1$, then the restriction of $\preceq$ to
$Stab_{A(\Gamma)}(F(0))$ is non-isolated, thus $\preceq$ itself is
non-isolated.

Now, if $Stab_{A(\Gamma)}(F(0))$ has rank exactly $1$, then the
order restricted to $Stab_{A(\Gamma)}(F(0))$ is completely
determined by the sign of any given non-trivial element, say
$A(t)\in Stab_{A(\Gamma)}(F(0))$. Assume $t\succ0$, then, because
$A(t)$ acts as a non-trivial homothety, there exists $x\in \R$ such
that $A(t)(x)>x$. It follows that $\preceq$ coincides with
$\preceq_{(x_1,x_2)}$, the ordering induced form the action $A$
where $x_1=F(0)$ and $x_2=x$. So $\preceq$ is non-isolated by
Proposition \ref{prop affine}. $\hfill\square$

\subsection{Conradian orderings}
\label{sec conrad}

There is a special type of left-ordering, introduced in
\cite{conrad}, which will be very important in our proof of Theorem
\ref{teo main}. These are the so called {\em Conradian} orderings,
which are left-orderings satisfying the following additional
property\footnote{In fact, in \cite{conrad} the required property is
$f\succ id, \; g\succ id \Rightarrow \exists n\geq 1$ such that
$fg^n\succ g$. The fact that $n=2$ is enough is a result from Navas
and Jim\'enez, see \cite{Navas}.}: $$f\succ id\;, \;\; g\succ id\;
\Rightarrow fg^2\succ g.$$

It turns out that Conradian orderings have a very interesting
dynamical counterpart. Recall from \cite{Navas}, that $f,g\in
Homeo_+(\R)$ are said to be {\em crossed}, if one of them, say $g$,
has a {\em domain} $I$ (that is, an open interval, not necessarily
bounded, which is fixed by $g$, and on which $g$ acts without fixed
points), such that $f(I)$ is not equal, nor disjoint to $I$. A group
$G\subset Homeo_+(\R)$ is said to act without crossings, if it does
not contains crossed elements.

\vsp

\begin{center}
\begin{tikzpicture}

\draw[<->] (-4,-1)-> (-1,2);

\draw [-] (-4,0.5) .. controls (-3,1) .. (-1.5,1.5); \draw [-]
(-4,-0.5)->(-1.5,2);

\draw (-4,0) node{$f$};\draw (-3,1.2) node{$g$};

\draw[<->] (1,-1)-> (4,2);

\draw [-] (1.5,-0.5) .. controls (2,1) .. (3.5,1.5); \draw [-]
(1.8,-1)->(3.5,1.5);

\draw (2.8,-0.2) node{$f$};\draw (1.5,0.3) node{$g$}; \draw[fill]
(3.5,1.5) circle (1.5 pt);

\draw[<->] (6,-1)-> (9,2);

\draw [-] (6.5,-0.5) .. controls (7,1) .. (8.5,1.5); \draw [-]
(6.8,-1)->(8.5,2);

\draw (8.1,1.7) node{$f$};\draw (6.5,0.3) node{$g$};

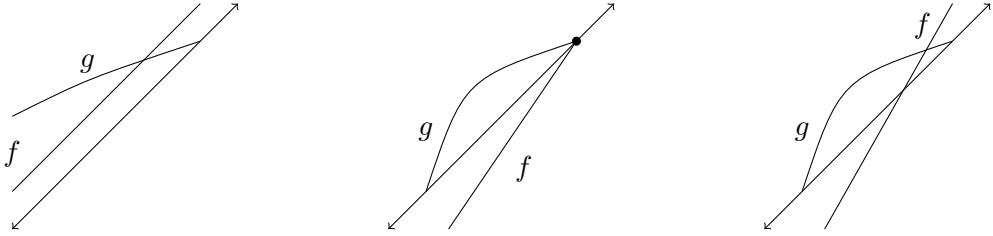
\draw (2.5,-1.7) node{{\small Figure 1: Three different crossings}};

\end{tikzpicture}
\end{center}

\begin{thm}[Navas \cite{Navas}] \label{teo navas}
{\em The dynamical realization of a Conradian ordering on $G$ is an
action without crossings. Conversely, an induced ordering from an
action without crossings is Conradian.}
\end{thm}

The above theorem implies rather easily that, in a dynamical
realization of a Conradian ordering, the set of elements having
fixed point is a subgroup (obviously normal). With this, together
with a theorem of H\"older stating that every group acting freely on
the real line is Abelian (see for instance \cite{ghys}), one can
deduce (see for instance \cite{Navas}. Compare with \cite{conrad}.)

\begin{cor}\label{coro conrad dyn real} {\em Let
$G$ be a countable group, $\preceq$ be a Conradian ordering of $G$,
and $N\subset G$ be the set of elements having a fixed point in the
dynamical realization of $(G,\preceq)$. Then $N$ is a normal
subgroup of $G$. Moreover, if there is $g\in G$ having no fixed
point (for instance if $G$ is finitely generated), then $G/N$ is a
non-trivial torsion-free Abelian group which acts freely on the
(non-empty) set of global fixed points of $N$. }
\end{cor}

Let us state here a last corollary which plays a crucial role in the
proof of Theorem \ref{teo main}. It describes some constrains on the
dynamics of a group acting on the real line, when it has a normal
subgroup acting without crossings. More precisely

\begin{cor}\label{cor no strong crossing}{\em
Let $(G,\preceq)$ be a left-ordered group, let $H$ be a normal
subgroup such that $(H,\preceq)$ is Conradian, and consider a
dynamical realization of $(G,\preceq)$. Let $f\in G$, $g\in H$ and
let $I$ be a minimal open interval fixed by $g$. Then one of the
following holds.
\begin{itemize}
\item $f(I)=I$ or
\item $f(I)$ is disjoint from $I$, or
\item (up to changing $f$ by its inverse)
$\overline{I}\subset f(I)$. In this last case we say that $f$ acts
as a {\em dilation} on $I$.
\end{itemize}}
\end{cor}

\noindent{\em Proof:} Notice that $f(I)$ is a domain of
$g^f=fgf^{-1}\in H$.  Hence, it follows from Theorem \ref{teo navas}
that $g$ and $g^f$ are not crossed. In particular, $I$ and its image
by $f$ are either disjoint or one is contained in the other. Indeed,
if this is not the case, then $f(I)$ would not be fixed, nor moved
disjointly by $g$. So, up to changing $f$ by $f^{-1}$, we can assume
that $I \subseteq f(I)$.

\vsp

We still have to rule out the possibility that these two intervals,
although different, share a common extremity. This is again easy.
Indeed, suppose it is the case that $I$ and $f(I)$ share a common
extremity but $I\subsetneq f(I)$. Then $g^f$ can not move $I$
disjoint from itself. But, on the other hand, since $g^f$ have no
fixed points inside $f(I)$, we have that $I$ can not be fixed by
$g^f$, contradicting Theorem \ref{teo navas}. $\hfill\square$

\section{The space of orderings of $SOL$} \label{sec examples}

\subsection{The space of orderings of $SOL$ is a Cantor set}
In this section we treat the simplest case of a non-virtually
nilpotent polycyclic group, namely  the group $SOL=\Z^2\rtimes_T
\Z$, where $T$ is an hyperbolic matrix (that is, a matrix in
$SL_2(\Z)$ having trace greater than $2$). In particular, $T$ has
two irrational eigenvalues. We begin by proving that this group has
no isolated orderings, and then we provide a quite explicit
description of its set of orderings (describing for instance its
bi-invariants orderings, which are left-orderings whose positive
cone is preserved under conjugation).

We denote by $H$ the derived subgroup of $SOL$, which is isomorphic
to $\Z^2$, and by $t$ the element of $\Z$ acting on $H$ as $T$. Let
$\preceq$ be a left-ordering on $SOL$, and consider its dynamical
realization.

Since $H$ is Abelian, the set of elements in $H$ acting with fixed
points form a subgroup $H'$. This subgroup, being finitely
generated, actually has a global fixed point. But because $T$ is
$\Q$-irreducible, this subgroup is either trivial or must have
finite index. In the latter situation, we have that $H'=H$ as every
global fixed point of $H'$ has a finite $H$-orbit, so must be fixed
by $H$. We therefore have two cases to consider, namely
\begin{itemize}

\item{\bf Case 1.} {\em $H$ has a global fixed point.}

\vsp

Let $I$ be the maximal open interval around $0$ without global fixed
point of $H$. Since $H$ is normal in $SOL$, and $SOL$ acts
without global fixed points, we have that the set of global fixed
points of $H$ is permuted by $SOL$ and therefore must be
infinite. In particular, $I$ is a bounded interval which is either
fixed or moved disjointly from itself by the action of $SOL$.
Therefore $H=Stab(I)$ is convex. Moreover, since $H$ has rank
two, the restriction of $\preceq$ to $H$ is non-isolated. Hence,
Proposition \ref{prop convex} implies that $\preceq$ itself is
non-isolated. Observe moreover that $\preceq$ is Conradian.

\item{\bf Case 2.} {\em $H$ has no global fixed point.}

\vsp

It follows that $H$ acts freely on the real line, and so, by
H\"older's Theorem \cite{ghys}, it is semi-conjugated to a group of
translation. Now recall that Lebesgue's measure is the unique
measure, up to scalar multiple, preserved by $\Z^2$ acting
faithfully by translations. In particular $H$ preserves a unique
atomless Radon\footnote{Recall that a {\em Radon} measure is a
measure giving finite mass to compact sets.} measure $\mu$.
Moreover, the hyperbolicity of $T$ implies that $t$ does not
preserves the measure $\mu$, but it acts on it as a dilation by one
of the two eigenvalues of $T$. We therefore obtain a faithful
embedding of $SOL$ in $Aff_+(\R)$ which is realized by a
semi-conjugation (see equation (\ref{eq affine action}) and (\ref{eq
semiconj to affine}) respectively). We then conclude from Corollary
\ref{cor:affine}. We observe that in this case, the proof of
Corollary \ref{cor:affine} shows that the ordering $\preceq$ is
realized as an induced ordering from the associated affine action.

\end{itemize}

\vs

\subsection{Description of the space of orderings.}
It follows from the previous analysis that there are two types of
orderings on $SOL$: those which are Conradian, and those which are
induced by affine actions.

\begin{itemize}
\item{\bf Conradian orderings.}
These always form a closed subset \cite{Navas}. Here, Conradian
orderings are exactly those for which the normal subgroup $H\simeq
\Z^2$ is convex. Therefore in $\mathcal{LO}(SOL)$ there are two
copies of the Cantor set $\mathcal{LO}(\Z^2)$, each of which
corresponding to a choice of sign for $t$. Let us briefly recall the
description of the space of left-orderings of $\Z^2$. First observe
that each oriented line passing through the origin in $\R^2$
delimits a unique half-plane (say the one on its right) defining the
positive cone of some pseudo-ordering on $\Z^2$. The set of elements
of $\Z^2$ which belong to the line form a cyclic convex subgroup,
which is trivial precisely when the pseudo-ordering is an ordering
(this happens exactly when the slope of the line is irrational). It
follows that this space of pseudo-orderings is naturally
parametrized by the unit circle. When the slope is rational, one
needs to specify an ordering of the convex subgroup, which is
determined by a sign. Therefore the space of orderings of $\Z^2$ can
be parametrized by a Cantor set obtained by ``doubling" each
rational point of the circle (see \cite{sikora}). For simplicity in
Figure 2, we ignore this ``blow up" procedure and represent each
copy of $\mathcal{LO}(\Z^2)$ as a ``vertical" circle.

\item{\bf Bi-orderings.}
Now observe that among these orderings on $\Z^2$, those which are
invariant under conjugation by $t$ are precisely those corresponding
to lines which are eigendirections of the matrix $T$. The
corresponding orderings of $SOL$ are those which are bi-invariant.
Taking into account the choices of orientations, this gives
precisely eight\footnote{A classification of finitely generated
groups admitting only finitely many bi-ordering can be found in
\cite{botto rhemtulla}.} bi-orderings.

\item{\bf Affine orderings.}
The complement of Conradian orderings, namely those induced by
affine actions, is an open subset with eight accumulation points,
namely $SOL$'s bi-orderings. We represent these affine orderings by
four copies of $\R$, compactified at $\pm \infty$ by pairs of
bi-orderings corresponding to two different eigendirections. This
requires some explanation.

First, these four intervals are to be thought as Cantor sets.
Similarly to our description of orderings of $\Z^2$, one can first
consider pseudo-orderings $\preceq_x^{(i)}$, $i=1\ldots 4$, induced
by the orbit of one point, $x$ under an affine action of $SOL$ on
the line. Such pseudo-ordering is determined by the following data:
first, choose an orientation of the line, then one needs to specify
the action of $t$ by multiplication by one of the two eigenvalues of
$T$. Hence these pseudo-orderings are naturally parametrized by four
copies of the real line. Note that $\preceq_x^{(i)}$  is an ordering
precisely when the stabilizer of $x$ is trivial, which happens on
the complement of some dense countable subset $D\subset \R$ (this
subset corresponds to the possible values of translations of
elements of $\Z^2$). Otherwise, to define an ordering on $SOL$, one
needs to specify an orientation on the convex cyclic subgroup of
homotheties fixing $x$. Therefore the subset of affine orderings is
locally a Cantor set, obtained by doubling points belonging to $D$
in each of the four copies of $\R$.

Consider one of these intervals, corresponding to an action where
$t$ acts by dilation (i.e.\ $t^{-1}$ acts by contraction), and let
$\preceq_x^{(1)}$ be the pseudo-ordering associated to the orbit of
$x$. Observe that when $x$ goes to $\infty$, the action of $t$
becomes predominant over translations, so that $\preceq_x^{(1)}$
converges to an ordering where translations form a convex subgroup.
Since the restriction of $\preceq_x^{(i)}$ to $H$ is conjugation
invariant, one easily checks that this limiting ordering is one of
our 8 bi-orderings.

\end{itemize}

Based on this description, it is not difficult to describe the
dynamics of the action of $SOL$ on its space of ordering, see for
instance Remark \ref{rem covariant actions}. We leave this as an
exercise to the reader.

\begin{center}
\begin{tikzpicture}

\draw[dashed] (-2,1) -- (2,1); \draw[dashed] (-2,0) -- (2,0);
\draw[dashed] (-2,-1) -- (2,-1); \draw[dashed] (-2,-2) -- (2,-2);

\draw[fill] (-2,1) circle (1.3 pt); \draw[fill] (2,1) circle (1.3
pt); \draw[fill] (-2,0) circle (1.3 pt); \draw[fill] (2,0) circle
(1.3 pt); \draw[fill] (-2,-1) circle (1.3 pt); \draw[fill] (2,-1)
circle (1.3 pt); \draw[fill] (-2,-2) circle (1.3 pt); \draw[fill]
(2,-2) circle (1.3 pt);

\draw[dashed] (-2,1) -- (-2,-2);\draw[dashed] (2,1) -- (2,-2);
\draw[dashed] (-2,1) .. controls (-3,1.5) and (-3,-2.5) .. (-2,-2);
\draw[dashed] (2,1) .. controls (3,1.5) and (3,-2.5) .. (2,-2);

\draw (0,-2.72) node { {\small Figure 2: The space of left-orderings
of $\Z^2\rtimes_T\Z$.}};

\end{tikzpicture}
\end{center}

\section{Plante's action of $\Z\wr\Z$.}\label{sec Plante}

The case of $SOL$ which has been described in the previous section
is a good representative of what happens for solvable groups with
finite rank. The goal of this subsection is to illustrate the
difficulties arising when dealing with  solvable groups of infinite
rank. For simplicity, we shall focus on Abelian-by-cyclic groups
$H\rtimes \Z$. The prototypical example of such group is
$\Gamma=\Z\wr\Z= \oplus_{i\in \Z}\Z \rtimes \Z$, where $\Z$ acts by
shifting the indices in $H=\oplus_{i\in \Z}\Z$.

In the finite rank case, we only have to consider two cases: either
$H$ has a global fixed point, or (has an element which) acts freely.
Here by contrast, since $H$ is not supposed to have finite rank, a
third case can happen for which all elements of $H$ have fixed
points, while $H$ does not have any global fixed point\footnote{In
the proof of Theorem \ref{teo main}, the first one of these three
cases (concerning Conradian orderings) is treated implicitely as an
initial step in our argument.}. In \cite[Sec. 5]{plante}, Plante
describes an explicit action of  $\Z\wr\Z$ corresponding to this
third case (in particular such an action does not quasi-preserve any
Radon measure).

We shall now recall the main properties of Plante's action of
$\Z\wr\Z$ which are responsible for the fact that the correspondings
induced orderings are non-isolated. Let $t$ be a generator of the
cyclic group acting on $H$. For this action, each non-trivial
element $h$ of $H$ preserves a compact interval containing $0$ in
its interior. Let $I_h$ be the minimal such interval. These
intervals are nested, and their interesection is reduced to $\{0\}$.
The element $t$ acts as an expanding homeomorphism having $0$ as its
unique fixed point, and such that $t(I_h)=I_{tht^{-1}}$. Somehow,
this action reminds of an affine action where $t$ would be some kind
of homothety, while elements of $H$  would play the role of
translations. It turns out that similarly to its affine cousin, this
ordering can be approximated by its conjugates. Indeed, take for
instance the pseudo-ordering induced by the orbit of $0$. The cyclic
subgroup generated by $t$ being the minimal convex subgroup, this
pseudo-ordering can be completed to an ordering by specifying the
sign of $t$. One easily checks that if $t\succ 0$ (resp.\ $t\prec
0$), then for any sequence of points $x_n$ converging to zero from
the right (resp.\ left), the sequence of orderings $\preceq_{x_n}$
converges to $\preceq$, while being distinct from it.

\begin{center}
\begin{tikzpicture}

\draw (1,2.7) node{\small $t$}; \draw[<->] (-2.5,-2.5)-> (2.5,2.5);
\draw[-] (-1.5,-3)-> (1.5,3);

\draw (-1.9,0) node{\small $h_0$}; \draw[-] (-2,-2) .. controls
(-1.34,1.34).. (2,2);\draw (2,2) .. controls (2.1,2.4) ..
(2.2,2.5);\draw (-2,-2) .. controls (-2.1,-2) .. (-2.5,-2.2);

\draw (-0.7,0.6) node{\small $h_1$}; \draw[-] (-1,-1) .. controls
(-0.6,0.6).. (1,1); \draw[-] (-1.76,-1.76) .. controls
(-1.65,-1.1).. (-1,-1); \draw[-] (1,1) .. controls (1.15,1.6)..
(1.76,1.76); \draw(1.76,1.76) .. controls (1.8, 1.9).. (1.92,1.92);
\draw(-1.92,-1.92) .. controls (-1.9, -1.8).. (-1.76,-1.76);

\draw[-] (-1.56,-1.56) .. controls (-1.5,-1.27).. (-1.2,-1.2);
\draw[-] (-0.8,-0.8) .. controls (-0.7,-0.5).. (-0.45,-0.45);
\draw[-] (-0.45,-0.45) .. controls (-0.3,0.3).. (0.45,0.45);\draw[-]
(0.45,0.45) .. controls (0.5,0.7).. (0.8,0.8); \draw[-] (1.2,1.2) ..
controls (1.3,1.45).. (1.56,1.56);

\draw[-] (-0.23,-0.23) .. controls (-0.1,0.1).. (0.23,0.23);

\draw[dashed](-1.76,-1.76)->(-1.76,-1)->(-1,-1)->(-1,1)->(1,1)->(1,1.76)->
(1.76,1.76); \draw[dashed](1,1)->(1,2)->(2,2);\draw[dashed]
(-2,-2)->(-1,-2)->(-1,-1);

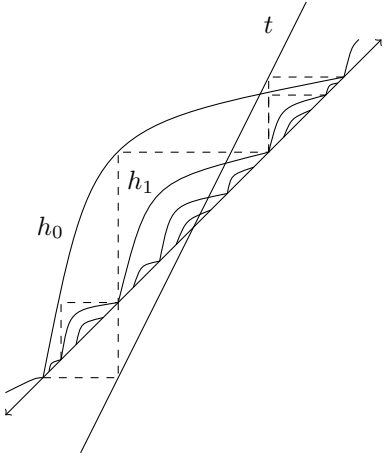
\draw (0,-3.5) node{\small Figure 3: Plante's action:
$h_1=th_0t^{-1}$.};

\end{tikzpicture}
\end{center}


\section{Proof of Theorem \ref{teo main}}

\label{sec Proof of Main thm}


Let $\preceq$ be a left-ordering on $\Gamma$. We assume that the
convex series of $(\Gamma,\preceq)$ is finite, because otherwise,
from Corollary \ref{cor infinite convex series}, $\preceq$ is
non-isolated. Say the convex series is
$$\{ id\}=C_n\subset C_{n-1}\subset \ldots \subset C_{0}=\Gamma.$$

\vsp

We let $T=C_{j+1}\subseteq \Gamma$ be the {\em Conradian soul} of
$(\Gamma,\preceq)$, which is the maximal convex subgroup on which
the restriction of $\preceq$ is Conradian. By \cite{Navas, Nav-Riv
crossings}, if $T$ is not a Tararin group, then $\preceq$ is
non-isolated. So we shall assume that $T$ is a Tararin group.

\vsp

If $T=\Gamma$, then we are done: $\Gamma$ admits only finitely many
left-orderings. So we suppose that $T=C_\ell$ is a proper convex
subgroup of $\Gamma$. We will show that the restriction of $\preceq$
to $C_{\ell-1}$ is non-isolated. Hence, in light of Proposition
\ref{prop convex}, there is no harm in assuming that there are no
convex subgroups between $T$ and $\Gamma$.

\vsp

To show that $\preceq$ is non-isolated, we consider a dynamical
realization of $(\Gamma,\preceq)$. We let $I_T$ be the minimal
interval stable by $T$ and containing $0$. The following is a direct
consequence of the fact that $T$ is convex.

\begin{lem}\label{lem:T}
{\em Every element of $\Gamma$ either fixes, or moves $I_T$
disjointly from itself. In particular, the stabilizer of $I_T$ is
exactly $T$.}
\end{lem}
\vsp

\vsp

We let $\widetilde{\Gamma}$ be a finite index,
normal, solvable subgroup of $\Gamma$. We let
$\widetilde{\Gamma}^0=\widetilde{\Gamma}$ and
$\widetilde{\Gamma}^i=[\widetilde{\Gamma}^{i-1},\widetilde{\Gamma}^{i-1}]$
be the associated derived series
$$\{ id\}=\widetilde{\Gamma}^k\lhd \widetilde{\Gamma}^{k-1}\lhd \ldots
\lhd \widetilde{\Gamma}^1\lhd \widetilde{\Gamma}\lhd \Gamma.$$ Note
that each $\widetilde{\Gamma}^i$ is normal in $\Gamma$. We fix once and
for all $i$, being the minimal index such that
$\widetilde{\Gamma}^{i}\subseteq T$. Since $T$ is a proper convex
subgroup, we have that $k\geq i \geq 1$. In a diagram

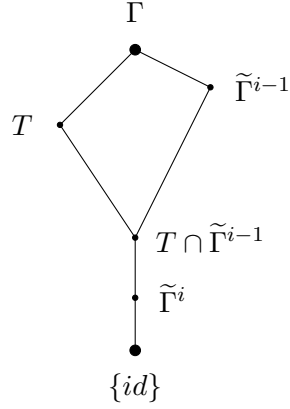
\begin{center}
\begin{tikzpicture}

\filldraw[black] (0,2) circle (2 pt);\filldraw[black] (0,-2) circle
(2 pt); \filldraw[black] (-1,1) circle (1 pt); \filldraw[black]
(1,1.5) circle (1 pt); \filldraw[black] (0,-0.5) circle (1
pt);\filldraw[black] (0,-1.3) circle (1 pt);

\draw[black] (0,2.5) node{$\Gamma$} ; \draw[black] (0,-2.5)
node{$\{id\}$} ; \draw[black] (0.5,-1.3)
node{$\widetilde{\Gamma}^{i}$} ; \draw[black] (-1.5,1) node{$T$} ;
\draw[black] (1.7,1.5) node{$\widetilde{\Gamma}^{i-1}$} ;
\draw[black] (1,-0.5) node{$T\cap\widetilde{\Gamma}^{i-1}$} ;

\draw [-](0,-2)-> (0,-0.5)->(-1,1)->(0,2)->(1,1.5)->(0,-0.5);

\end{tikzpicture}
\end{center}

\begin{lem}\label{lem:the restriction is Conrad}{\em The order
restricted to $\widetilde{\Gamma}^{i-1}$ is Conradian.}
\end{lem}

\noindent{\em Proof:} Indeed, $\widetilde{\Gamma}^{i-1}\cap T$ is
convex and normal in $\widetilde{\Gamma}^{i-1}$. Moreover, its
quotient is Abelian, so it admits only Conradian orderings. The
lemma then follows from Remark \ref{rem:ConradExtension}.
$\hfill\square$

\vsp

\begin{lem}\label{lem:the action is cofinal}{\em The orbit of
$0$ under $\widetilde{\Gamma}^{i-1}$ accumulates on $\pm\infty$.}
\end{lem}

\noindent{\em Proof:} Let $I$ be the smallest open interval
containing $0$, and stable under $\widetilde{\Gamma}^{i-1}$. Since
$\widetilde{\Gamma}^{i-1}$ is normal in $\Gamma$, $I$ is either
fixed or moved disjointly by any $\gamma\in \Gamma$. In particular,
$Stab_\Gamma(I)$, the stabilizer of $I$, is a convex subgroup. Now,
if $I\not= \R$, then the maximality of $T$ implies
$Stab_\Gamma(I)\subseteq T$. But this implies
$\widetilde{\Gamma}^{i-1}\subset T$, which is contrary to our
assumptions. $\hfill\square$

\vs

We have two cases to analyze in order to prove Theorem \ref{teo
main}.

\begin{itemize}

\item{\bf Case 1.}
{\em Suppose there is $g_0\in
\widetilde{\Gamma}^{i-1}$ having no fixed points.} Such case occurs for instance if
$\widetilde{\Gamma}^{i-1}$ has finite rank. Combining the results from \cite{Navas,plante} we
obtain the following proposition (by measure we shall implicitely
mean a Radon measure which is finite on compact sets).

\begin{prop}\label{prop semi conj to affine}{\em The
action of $\Gamma$ on the real line is semi-conjugated to a
non-Abelian affine action $\varphi:\Gamma\to Aff_+(\R)$.}
\end{prop}

\noindent{\em Proof:} Since the action of $\widetilde{\Gamma}^{i-1}$
is Conradian and there is $g_0\in \widetilde{\Gamma}^{i-1}$ without
fixed points, Corollary \ref{coro conrad dyn real} ensures that
$\widetilde{\Gamma}^{i-1}$ has a maximal proper convex subgroup
$N\supseteq \widetilde{\Gamma}^i$, which is normal. In particular,
$N$ fixes some open bounded interval $I_N$ around $0$. The action of
$\widetilde{\Gamma}^{i-1}$ is semi-conjugated to an action factoring
through $A=\widetilde{\Gamma}^{i-1}/N$. Since $A$ is Abelian and has
an element acting without fixed point, \cite[Proposition
3.1]{plante} implies that the corresponding action has an invariant
measure. Lifting back this measure yields a measure $\mu$ which is
preserved by $\widetilde{\Gamma}^{i-1}$ in the original action.

Since $g_0$ acts without fixed points, the {\em translation number
homomorphism} $\tau_\mu:\widetilde{\Gamma}^{i-1}\to \R$, given by
$\tau_\mu(g)=\mu((0,g(0)])$, is non trivial (here and below, we use
the convention $\mu([x,y])=-\mu([y,x])$ for $y<x$). It now follows
from \cite[Lemma 4.2 and 4.3]{plante} that there is a measure, which
we still call $\mu$, which is {\em quasi-preserved} by $\Gamma$,
meaning that for every $\gamma\in \Gamma$, there is a positive real
number $\lambda_\gamma$ such that $\gamma_*(\mu)=\lambda_\gamma \mu$
(where $\gamma_*(\mu)(X):=\mu(\gamma^{-1}(X))$, $X\subseteq \R$). In
this way we have a homomorphism $\varphi:\Gamma\to Aff_+(\R)$, which
extends $\tau_\mu$, given by
\begin{equation}\label{eq affine action}\varphi(\gamma)(x)=\frac{1}{\lambda_\gamma} \, x +
\mu((0,\gamma(0)]).\end{equation} This affine action is
semi-conjugated to the original dynamical realization action of
$(\Gamma,\preceq)$. Indeed, if for $x\in \R$ we let
$F(x)=\mu((0,x])$, then a direct computation shows that
\begin{equation}\label{eq semiconj to affine} F(\gamma(x))=
\varphi(\gamma)(F(x)).\end{equation} It only remains to check that
$\varphi(\Gamma)$ is non-Abelian. To this end we let $I_\mu:=(a,b)$,
where $a=\sup\{x<0\mid x\in supp(\mu)\}$, and $b=\inf=\{x>0\mid x\in
supp(\mu)\}$. Then, since $\mu$ is quasi-preserved (so in particular
its support $supp(\mu)$ is preserved) we have that
$Stab_\Gamma(I_\mu)$ is either fixed or moved disjointly. In
particular, $Stab_\Gamma(I_\mu)$ is a proper convex subgroup. So,
$Ker(\varphi)\subseteq Stab_\Gamma(I_\mu)\subseteq T$. Therefore if
the affine action of $\Gamma$ was Abelian, then $T$ would be a
normal and co-Abelian subgroup, so Remark \ref{rem:ConradExtension}
would imply that $\preceq$ is Conradian, which is contrary to our
assumptions. $\hfill\square$

\vsp

\begin{lem}\label{lem kernel is convex} {\em The kernel of $\varphi$ is a
convex subgroup of $(\Gamma,\preceq)$.}
\end{lem}

\noindent{\em Proof:} We keep the notations of the proof of
Proposition \ref{prop semi conj to affine}.

We first claim that $Stab_\Gamma(I_\mu)=T$ (equivalently
$I_\mu=I_T$). Indeed, let $\preceq'$ be the pseudo-ordering of
$\Gamma$ induced by $\varphi$ at $F(0)=F(I_\mu)$. Since $\varphi$ is
an affine action, $\preceq'$ has only one convex subgroup, namely
$\{\gamma\in \Gamma\mid \varphi(\gamma)(F(0))=F(0)\}
=Stab_\Gamma(I_\mu)$. However, equation (\ref{eq semiconj to
affine}) implies that $\preceq'$ is the quotient of $\preceq$ by
$Stab_\Gamma(I_\mu)$: $\varphi(g)\succ' id \Rightarrow g\succ id$.
Thus, the claim follows from Remark \ref{rem convex}.

We now show that $Ker(\varphi)$ is convex in $(\Gamma,\preceq)$.
First observe that the previous claim implies that $\varphi(T)$ is
the Abelian subgroup of homotheties centered at $F(0)$. If it was
trivial, then we would have that $T=\ker\varphi$ is convex. Let us
therefore suppose that $\varphi(T)$ is non-trivial.

We let $\{id\}=T_0\lhd T_1\lhd \ldots \lhd T_m=T$ be the convex
series of the Tararin group $T$. Recall that $T_i/T_{i-1}$ has rank
$1$ and that the action of $T_{i+1}$ on $T_i/T_{i-1}$ is by
multiplication by some negative number. In particular $T/T_{m-1}$ is
the unique torsion-free Abelian quotient of $T$. It follows that
$T/T_{m-1}=\varphi(T)$ and $\ker \varphi=T_{m-1}$ is convex.
$\hfill\square$

Now the proof of Theorem \ref{teo main} in case 1 follows from
Corollary \ref{cor:affine}.

\item{\bf Case 2.} {\em Suppose every element $g\in
\widetilde{\Gamma}^{i-1}$ has a bounded {\em domain} ({\em i.e.} an
open interval fixed by $g$ with no fixed point of $g$ in its
interior) $I_g$ around $0$.}

Our proof of Theorem \ref{teo main} in Case 2 consists in showing
that $\preceq$ can be approximated by a left-ordering induced from
the dynamical realization of $(\Gamma,\preceq)$, where the first
reference point is chosen outside but very close to $I_T$. For this
purpose, we shall prove that the action is quite similar to the one
described in \S \ref{sec Plante} (except that there, $I_T$ was
reduced to a point).

 It follows from Lemmas \ref{lem:the restriction is Conrad} and
\ref{lem:the action is cofinal} that the union of the $I_g$ is all
of $\R$. In the sequel, we exploit the facts that
$\widetilde{\Gamma}^{i-1}$ is normal in $\Gamma$ and that the order
restricted to $\widetilde{\Gamma}^{i-1}$ is Conradian to give a more
precise picture of the action.  First, Corollary \ref{cor no strong
crossing} immediately implies

\begin{lem} \label{lem no strong crossings}{\em Let $f\in \Gamma$ and
$g\in \widetilde{\Gamma}^{i-1}$. Then one, and only one of the
following happen \begin{itemize} \item $f(I_g)=I_g$ or
\item $f(I_g)$ is disjoint from $I_g$ or either
\item (up to changing $f$ by its inverse)
$\overline{I_g}\subset f(I_g)$, i.e.\ $f$ acts as dilation on $I_g$.\end{itemize}}

\end{lem}

\vs

The following corollaries are easy and left to the reader.

\begin{cor}\label{cor:T either dilates or fixes} {\em Let $g\in
\widetilde{\Gamma}^{i-1}\setminus T$, and $I_g$ its domain. If $t\in
T$, then $t$ either fixes $I_g$ or acts on it as a dilation.}
\end{cor}

\begin{cor}\label{cor: union/intersection} {\em If $I$ is an
interval obtained as a union or an intersection of $I_g$'s for $g\in
\widetilde{\Gamma}^{i-1}$, then a weak form of Lemma \ref{lem no
strong crossings}  holds for $I$. Namely every element $f\in \Gamma$
either moves $I$ disjointly from itself, or up to replacing $f$ by
its inverse, $I \subset f(I)$ (if the intersection is strict, we say
that $f$ weakly dilates $I$).}
\end{cor}

Using the fact that Tararin
groups have finite rank, we now obtain a useful description of $I_T$..

\begin{lem}\label{lem:intersection}
{\em The intersection of  $I_g$ for $g\in \widetilde{\Gamma}^{i-1}
\setminus T$ coincides with $I_T$.}
\end{lem}

\noindent {\em Proof:} First observe that Theorem \ref{teo navas}
implies that the $I_g$'s, for $g\in \widetilde{\Gamma}^{i-1}\setminus
T$  are totally ordered for the inclusion.  On the other hand, by Corollary \ref{cor: union/intersection}, if no element weakly dilates $J$, then the stabilizer of
$J$ is convex. Thus it contains $T$ by Corollary \ref{cor:T either
dilates or fixes} (which also applies to $J$), so it must be equal to $T$. But this implies that
$J=I_T$.

It is therefore enough to prove that no element weakly dilates $J$.
Suppose by contradiction that $f^{-1}$ weakly dilates $J$ (i.e. $f$
weakly contracts $J$). Then there exists $g\in
\widetilde{\Gamma}^{i-1}\setminus T$ such that $J\cap I_{g^f}$ is
strictly contained in $J$, hence that $g^f$, and more generally
$g_n=g^{f^n}$ belongs to $T$ for all $n\geq 1$. Moreover, $f^n(I_g)$
is a minimal interval (not necessarily containing $0$) fixed by the
element $g_n$. Since $f$ acts as a dilation on $I_g$, we deduce that
$f^{n+1}(I_g)\subsetneq f^n(I_g)$. In particular  $g_k\in S_n\;
\forall k\geq n$, where $S_n=Stab_T(f^n(I_g))$. Let $x$ be a point
in the decreasing intersection of compact intervals
$\overline{f^n(I_g)}$. We have that $(S_n)_{n\geq 1}$ is a strictly
decreasing sequence of convex subgroups of $T$ for the
pseudo-ordering $\preceq_x$, violating the fact that $T$ has only
finitely many orderings (and more generally one could check that
this is incompatible with the fact that $T$ has finite rank).
$\hfill\square$

Since there are no proper convex subgroups above $T$, we have

\begin{lem}\label{lem:existence of a dilation} {\em For any
$g\in \widetilde{\Gamma}^{i-1}\setminus T$, there exists $f\in
\Gamma$ acting as a dilation on $I_g$. Moreover, for every $g, g'\in
\widetilde{\Gamma}^{i-1}\setminus T$ satisfying $I_g\subseteq
I_{g'}$, there exists $f\in \Gamma$ such that $I_{g'}\subseteq
f(I_g)$.}
\end{lem}

\noindent{\em Proof:} Looking for a contradiction, suppose that
there exists $g\in \widetilde{\Gamma}^{i-1}\setminus T$ such that
elements of $\Gamma$ either stablize $I_g$ or move it disjointly
from itself. It follows that the stabilizer $S$ of $I_g$ is convex.
By Corollary \ref{cor:T either dilates or fixes}, $S$ contains  $T$,
and the inclusion is strict as $S$ contains $g$. On the other hand
$S$ cannot be all of $\Gamma$ since $I_g$ is bounded: a
contradiction. This shows the first part of the lemma.

Let us prove the second statement of the lemma. Let
$\mathcal{F}=\{f\in \Gamma\mid f$ dilates $I_g\}$. First note that
for all $f\in \mathcal{F}$, we have $f(I_g)=I_{g^f}$. Therefore the
set of all $f(I_g)$ for $f\in \mathcal{F}$ is nested, let $I$ be its
union. We claim that $I$ is either fixed or moved disjointly.

Indeed, let us first suppose that there exists $h\in
\widetilde{\Gamma}^{i-1}$ such that $h(I)$ strictly contains $I$. A
continuity argument implies that there exits $f\in \mathcal{F}$ such
that $h\circ f(I_g)$ contains $f(I_g)$ (hence $I_g$) and is not
contained in $I$. The first of these statements implies that $h\circ
f$ belongs to $\mathcal{F}$, while the other implies that it does
not (by definition of $I$), so this case cannot occur. Hence, from
Corollary \ref{cor: union/intersection}, our claim follows. In
particular, the stabilizer of $I$ is a convex subgroup containing
$T$, so $I$ must be all of $\R$. Hence the lemma. $\hfill\square$

We have given a combinatorial description of the dynamics of
$\Gamma$ around $I_T$. We now exploit this description to
approximate $\preceq$. First we show that orderings induced by
points outside of $I_T$ are distinct from $\preceq$.

\begin{lem}\label{lem:distinct}
{\em For all $x$ not in $\overline{I_T}$, there exists $f\in\Gamma$
such that $f\prec_x 0$ and $f\succ 0$. In particular, any
left-ordering induced from the dynamical realization of
$(\Gamma,\preceq)$ with $x$ as first reference point is different
from $\preceq$.}
\end{lem}

\noindent{\em Proof:} Let $x\notin \overline{I_T}$, say on its left
 (the other case is symmetric). It results from Lemma
\ref{lem:intersection} that one can find $g,g'\in
\widetilde{\Gamma}^{i-1}\setminus T$ such that $I_g\subset I_{g'}$
and $x\in I_{g'}\setminus I_g$. On the other hand, Lemma
\ref{lem:existence of a dilation} provides us with an element $f$
such that $f(I_g)$ contains $I_{g'}$. In particular, if $f\succ id$,
then we are done because $f\prec_x id$. So we assume that $f\prec
id$. Up to changing $g$ by $g^{-1}$, we can assume that $g\succ id$.
In particular, $g(I_T)$ is moved to the right of $I_T$. This
implies, that for $n$ large enough, $g^{-n}fg^n \succ id$. But, in
the same time, $g^{-n}fg^n(I_g)=g^{-n} f(I_g)=g^{-n}
(I_{g^f})=f(I_g)$, where the last equality follows because $g$ and
$g^f$ are not crossed and $I_g\subset I_{g^f}$. Hence the lemma.
$\hfill\square$

\vsp

The following step consists in showing that when $x_n$ converges to
an end point of the interval $I_T$, then $\preceq_{x_n}$ converges
to $\preceq$ outside of $T$.

\begin{lem}\label{prop:approx}{\em If $x_n$ converges to an end point
$x$ of $I_T$, then $\preceq_{x_n}$ converges to $\preceq$ outside of
$T$. More precisely, for every $g\in\Gamma\setminus T$, we have that
$g\succ 0$ if and only if $g\succ_{x_n}0$ for all $n$ large enough.}
\end{lem}

\noindent{\em Proof:} Since the stabilizer of $I_T$ is precisely
$T$, given any $\preceq$-positive $g\in \Gamma\setminus T$, we have
that $g$ moves $I_T$ disjointly to its right. Therefore $g\succ 0$
if and only if $g\succ_{x}0$. The lemma then follows by continuity.
$\hfill\square$

\vsp

To prove Theorem \ref{teo main}, we are left to proving that for a
well-chosen sequence $(x_n)$ converging to an end point of $I_T$,
the orderings $\preceq_{x_n}$ converges to $\preceq$ {\it in
restriction to $T$}. We shall use in a crucial way the fact that
Tararin groups have ``very few" actions on the real line.  More
precisely, we have two following lemmas where the action considered
is still the dynamical realization of $(\Gamma,\preceq)$. Recall
that the convex series of $T$ is always given by $\{id\}=T_0\lhd
T_1\lhd \ldots\lhd T_{m-1}\lhd T_m=T$. We let $\gamma_T\in T$ be a
non trivial element in $T/T_{m-1}$ which acts on $T_{m-1}/T_{m-2}$
by a multiplication by a negative (rational) number.

\begin{lem}\label{lem:TararinFix}
{\em If $\gamma_T$ fixes some point $x$, then so does $T$. }
\end{lem}

\noindent{\em Proof:} Suppose $x$ is fixed by $\gamma_T$ but not by
$T$. We can then induce a left-ordering on $T$ with reference point
$(x,x_2,\ldots)$. In this ordering the stabilizer of $x$ is a proper
convex subgroup, which is impossible since it does not coincide with
any of the subgroups $T_i$. $\hfill\square$

\begin{lem}\label{lem:TararinOrder} {\em For any $y\in \R$, which is not
fixed by $\gamma_T$, there is $x$ between $\gamma_T^{-2}(y)$ and
$\gamma_T^2(y)$ such that $\preceq$ and $\preceq_{x}$ coincide over
$T_{m-1}$.}
\end{lem}

\noindent{\em Proof:} Recall that in a dynamical realization, the
set of fixed points of a non-trivial element has empty interior.
Since $T$ is countable, there is a point $z$ between $y$ and
$\gamma_T^{-1}(y)$  whose orbit under $T$ is free. In particular
$\preceq_z$ is a total ordering of $T$. Since $T_{m-1}$ is convex in
$\preceq_z$, there is a minimal open interval $I$ containing $z$ and
which is stabilized by $T_{m-1}$. Being moved disjointly from itself
by any non-trivial power of $\gamma_T$, $I$ contains at most one
point of the orbit of $y$ under $\langle \gamma_T\rangle.$ In
particular, $I$ is (strictly) contained between $\gamma_T^{-2}(y)$
and $\gamma_T(y)$. Now by Proposition \ref{prop two orbits of
Tararin}, there exists an element $g$ either in $T_{m-1}$ or in
$\gamma_TT_{m-1}$ such that $\preceq$  and
$g(\preceq_z)=\preceq_{gz}$ coincide over $T_{m-1}$. Clearly $x=gz$
satisfies the conclusion of the lemma.
 $\hfill\square$

\vs

The last step in the proof of Theorem \ref{teo main} is achieved by the following lemma.

\begin{lem} \label{lem approx over T}{\em There exists a sequence
$(x_n)$ converging to an end point of $I_T$ from outside such that
the induced left-ordering $\preceq_{x_n,\, 0}$ coincides with
$\preceq$ over $T$ for all $n$.}
\end{lem}

Indeed, combining Lemmas \ref{prop:approx} and \ref{lem approx over
T}, We have that $\preceq_{x_n,\, 0}$ converges to $\preceq$. On the
other hand, Lemma \ref{lem:distinct} shows that $\preceq_{x_n,\, 0}$
and $\preceq$ are different. This shows Theorem \ref{teo main} in
Case 2.

\vs

\noindent{\em Proof of Lemma \ref{lem approx over T}:} The idea is
to take $x\notin \overline{I_T}$, $x$ close to $\partial I_T$, such
that the sign of $\gamma_T$ is preserved and then to use  Lemmas
\ref{lem:TararinFix} and \ref{lem:TararinOrder}.

By Lemma \ref{lem:intersection}, there exists a sequence $g_n\in
\widetilde{\Gamma}^{i-1}\setminus T$ such that  $I_{g_n}$ converges to $I_T$.
\begin{itemize}
\item{\bf Subcase 1.} {\em  $\gamma_T$ dilates
$I_{g_n}$.}

Up to taking a subsequence, we can assume that $I_{g_n}\subset
I_{g_{n-1}}$. For concreteness we suppose $\gamma_T\prec id$, and
let $y_n$ be the left end-point of $I_{g_n}$, so that the sequence
$(y_n)$ converges to the left end-point $z$ of $I_T$ (in the case
$\gamma_T\succ id$, we take $y_n$ being right end-points). Since
$\gamma_T$ dilates $I_{g_n}$, $\gamma_T\prec_{y_n} id$. Let $x_n$ be
the sequence of points obtained from Lemma \ref{lem:TararinOrder}
between $\gamma_T^{-2}(y_n)$ and $\gamma_T^{2}(y_n)$ such that
$\preceq$ and $\preceq_{x_n}$ coincide in restriction to $T_{m-1}$.
By continuity of $\gamma_T$ and its inverse, the sequence $(x_n)$
converges to the same limit $z$. This shows the lemma in Subcase 1.

\item{\bf Subcase 2.} {\em  $\gamma_T$ (hence $T$) fixes $I_{g_n}$.}

\vsp

Again we let $x_n$ be the left end-point of $I_{g_n}$. Observe that
$\preceq_{x_n}$ is a pseudo-ordering on $\Gamma$, for which the
stabilizer $S_n$ of $x_n$ is obviously convex and, from Lemma
\ref{lem:TararinFix}, contains $T$. In particular, the left-ordering
$\preceq_{x_n,\, 0}$ coincides with $\preceq$ over $T$. This ends
the proof of Lemma \ref{lem approx over T}. $\hfill\square$
\end{itemize}
\end{itemize}


\begin{small}

\vspace{0.1cm}


\vspace{0.37cm}

\noindent Crist\'obal Rivas\\
\noindent UMPA, ENS-Lyon\\
\noindent Email: cristobal.rivas@ens-lyon.fr\\

\noindent Romain Tessra\\
\noindent UMPA, ENS-Lyon\\
\noindent Email: Romain.Tessera@ens-lyon.fr\\

\end{small}


\end{document}